\documentclass[11pt,oneside]{amsart}

\usepackage[noauthor,nofancy,nothm]{evan}

\usepackage[T1]{fontenc}
\usepackage{lmodern}

\usepackage[letterpaper,margin=1in,includehead,includefoot]{geometry}

\usepackage{setspace}
\usepackage[section]{placeins}

\usepackage{fancyhdr}
\fancypagestyle{stsbody}{%
  \fancyhf{}%
  \fancyhead[C]{\footnotesize\scshape \papertitle}%
  \fancyfoot[R]{\thepage}%
}
\usepackage{thmtools}   % amsthm itself is already built into amsart
\usepackage{graphicx}
\usepackage{tikz}
\usepackage{booktabs,graphicx}
\usepackage{pgfplots}
\usepackage{pgfplotstable}
\usetikzlibrary{arrows.meta,calc,shadings}
\pgfplotsset{compat=1.18}
\usetikzlibrary{angles,quotes,calc,patterns,decorations.pathreplacing}
\hypersetup{hypertexnames=false}

\declaretheorem[style=plain,name=Theorem,numberwithin=section]{theorem}
\declaretheorem[style=plain,name=Lemma,sibling=theorem]{lemma}
\declaretheorem[style=plain,name=Proposition,sibling=theorem]{proposition}

\declaretheorem[style=plain,name=Theorem,numbered=no]{theorem*}
\declaretheorem[style=plain,name=Lemma,numbered=no]{lemma*}
\declaretheorem[style=plain,name=Proposition,numbered=no]{proposition*}
\declaretheorem[style=plain,name=Corollary,numbered=no]{corollary*}
\declaretheorem[style=plain,name=Claim,numbered=no]{claim*}
\declaretheorem[style=plain,name=Conjecture,numbered=no]{conjecture*}
\declaretheorem[style=plain,name=Fact,numbered=no]{fact*}
\declaretheorem[style=plain,name=Assumption,numbered=no]{assume*}

\declaretheorem[style=definition,name=Definition,sibling=theorem]{definition}
\declaretheorem[style=definition,name=Example,sibling=theorem]{example}

\declaretheorem[style=definition,name=Definition,numbered=no]{definition*}
\declaretheorem[style=definition,name=Example,numbered=no]{example*}
\declaretheorem[style=definition,name=Algorithm,numbered=no]{algorithm*}
\declaretheorem[style=definition,name=Answer,numbered=no]{answer*}
\declaretheorem[style=definition,name=Exercise,numbered=no]{exercise*}
\declaretheorem[style=definition,name=Problem,numbered=no]{problem*}
\declaretheorem[style=definition,name=Question,numbered=no]{ques*}

\declaretheorem[style=remark,name=Remark,numbered=no]{remark*}
\declaretheorem[style=remark,name=Note,numbered=no]{note*}

\Crefname{assume}{Assumption}{Assumptions}
\Crefname{ques}{Question}{Questions}

\renewcommand{\vocab}[1]{\emph{#1}}

\DeclareMathOperator{\Opt}{Opt}

\newcommand{\fracpart}[1]{\left\langle #1 \right\rangle}

\makeatletter
\let\@setaddresses\relax
\makeatother

\newcommand{\credit}[1]{%
  \par\vspace{0.35\baselineskip}%
  \begin{spacing}{1}%
    \centering
    \begin{minipage}{\dimexpr\linewidth-2\captionindent\relax}%
      \centering\footnotesize\itshape #1%
    \end{minipage}%
    \par
  \end{spacing}}

\newcommand{\papertitle}{An Asymptotic for E. P. Stein Numbers}
\newcommand{\paperauthor}{Cameron Rampell}
\newcommand{\paperschool}{The Harker School}
\newcommand{\paperemail}{cameronrampell@gmail.com}
\providecommand{\paperemail}{}
\newcommand{\ifpaperemail}[1]{%
  \if\relax\detokenize{\paperemail}\relax\else #1\fi}

\title{\papertitle}
\author{\paperauthor}
\address{\paperschool}
\ifpaperemail{\email{\paperemail}}

\date{\today}


\begin{document}

%%%%%%%%%%%%%%%%%%%%%%%%%%%%%%%%%%%%%%%%%%%%%%%%%%%%%%%%%%%%%%%%%%%%%%%%%%%%%%%
%% TITLE PAGE -- separate, unnumbered, and not counted toward the 20-page
%% limit.  Set single-spaced: the 1.5-spacing rule applies to the body text.
%%%%%%%%%%%%%%%%%%%%%%%%%%%%%%%%%%%%%%%%%%%%%%%%%%%%%%%%%%%%%%%%%%%%%%%%%%%%%%%
\pagenumbering{gobble}
\thispagestyle{empty}
\begin{spacing}{1}
\begin{center}
  \vspace*{2in}

  {\LARGE\bfseries \papertitle \par}

  \vspace{2.5em}
  {\large \paperauthor \par}

  \vspace{0.6em}
  {\large \paperschool \par}

  \ifpaperemail{%
    \vspace{0.6em}
    {\normalsize \texttt{\paperemail} \par}%
  }

  \vfill
  {\today \par}
  \vspace{1in}
\end{center}
\end{spacing}
\clearpage

%%%%%%%%%%%%%%%%%%%%%%%%%%%%%%%%%%%%%%%%%%%%%%%%%%%%%%%%%%%%%%%%%%%%%%%%%%%%%%%
%% ABSTRACT PAGE -- separate, unnumbered, and not counted toward the 20-page
%% limit.  Typeset by hand because amsart's abstract environment prints only
%% as part of \maketitle, which the separate title page replaces.
%%%%%%%%%%%%%%%%%%%%%%%%%%%%%%%%%%%%%%%%%%%%%%%%%%%%%%%%%%%%%%%%%%%%%%%%%%%%%%%
\thispagestyle{empty}
\begin{center}
  {\large\bfseries Abstract\par}
\end{center}
\vspace{1.5em}

\textit{Epstein's put-or-take-a-square game} is an impartial combinatorial game
in which a player may, in each turn, add or subtract the largest square not
exceeding the current integer. We study \textit{E. P. Stein} numbers, the
positions from which both legal moves lead directly to positive squares, and
prove the following asymptotic formula for the counting function $S_2(X)$ of
these positions:
\[
  S_2(X)=\frac{1}{\pi^2}X^{1/4}\log X+O(X^{1/4}).
\]
% local tightness of asymptotic, optional

\vfill
\begin{spacing}{1}

\end{spacing}
\clearpage

%%%%%%%%%%%%%%%%%%%%%%%%%%%%%%%%%%%%%%%%%%%%%%%%%%%%%%%%%%%%%%%%%%%%%%%%%%%%%%%
%% CONTENT PAGES -- numbering restarts here, so the printed page numbers are
%% exactly the pages that count against the 20-page limit.
%%%%%%%%%%%%%%%%%%%%%%%%%%%%%%%%%%%%%%%%%%%%%%%%%%%%%%%%%%%%%%%%%%%%%%%%%%%%%%%
\pagenumbering{arabic}
\pagestyle{stsbody}

\section{Introduction and Motivation}

\textit{Nim} is a mathematical game of strategy in which two players take turns removing
objects from distinct heaps or piles; the player who removes the last object
wins. Game theorist Richard A. Epstein proposed a variant of Nim known as the
\textit{put-or-take-a-square game} \cite{winning_ways}, which we study in this paper.

The game starts with a positive integer. Each turn, writing $n$ for the
number the player to move is faced with, that player has two moves:
\begin{enumerate}
  \item Subtract $\floor{\sqrt{n}}^2$ from $n$; that is, subtract off $m^2$,
    where $m$ is the largest integer with $m^2 \leq n < (m + 1)^2$. We write
    $D(n)$ for this move.
  \item Add $\floor{\sqrt{n}}^2$ to $n$. We write $U(n)$ for this move.
\end{enumerate}
Formally, the two functions, for ``down'' and ``up'' respectively, are
\begin{align*}
  D(n) &= n - \floor{\sqrt{n}}^2, \\
  U(n) &= n + \floor{\sqrt{n}}^2.
\end{align*}
So from the position $n$, the set of positions reachable in one move is
$\Opt(n) = \{D(n), U(n)\}$. The game ends when a player moves to $0$, and that
player wins.

\begin{example}
  We illustrate by playing a put-or-take-a-square game starting with $n = 10$.
  \begin{enumerate}
    \item Player 1 starts with $n = 10$. Since $D(10) = 1$ and $U(10) = 19$,
      the moves are $\Opt(10) = \{1, 19\}$; say Player 1 chooses
      $U(10) = 19$.
    \item Player 2 now has $n = 19$, whose moves are
      $\Opt(19) = \{3, 35\}$; say Player 2 chooses $D(19) = 3$.
    \item Player 1 now has $n = 3$, whose moves are $\Opt(3) = \{2, 4\}$;
      say Player 1 chooses $U(3) = 4$.
    \item Player 2 now has $n = 4$, whose moves are $\Opt(4) = \{0, 8\}$.
      Player 2 chooses $D(4) = 0$ and wins the game.
  \end{enumerate}
\end{example}

\begin{figure}[ht]
\centering
\begin{tikzpicture}[
    >=Latex,
    font=\small,
    every node/.style={font=\small},
    mainnode/.style={
        circle,
        draw=blue!75!black,
        line width=1.0pt,
        fill=blue!6,
        minimum size=11.5mm,
        inner sep=0pt
    },
    sidenode/.style={
        circle,
        draw=gray!70,
        line width=0.9pt,
        fill=gray!6,
        minimum size=10.5mm,
        inner sep=0pt
    },
    mainedge/.style={
        ->,
        draw=blue!75!black,
        line width=1.1pt
    },
    sideedge/.style={
        ->,
        draw=gray!60,
        line width=0.85pt
    },
    playerlabel/.style={
        text=blue!75!black,
        font=\small
    }
]

% --- main nodes ---
\node[mainnode] (n10) at (0,0) {$10$};
\node[mainnode] (n19) at (2.3,-1.8) {$19$};
\node[mainnode] (n3)  at (0.6,-3.5) {$3$};
\node[mainnode] (n4)  at (2.4,-5.3) {$4$};
\node[mainnode] (n0)  at (1.0,-7.1) {$0$};

% --- side nodes ---
\node[sidenode] (n1)  at (-2.0,-1.8) {$1$};
\node[sidenode] (n35) at (4.2,-3.5) {$35$};
\node[sidenode] (n2)  at (-1.2,-5.3) {$2$};
\node[sidenode] (n8)  at (4.3,-7.1) {$8$};

% --- player labels ---
\node[playerlabel, above=2.5mm of n10] {Player 1};
\node[playerlabel, right=3.5mm of n19] {Player 2};
\node[playerlabel, right=3.5mm of n3]  {Player 1};
\node[playerlabel, right=3.5mm of n4]  {Player 2};
\node[playerlabel, below=1.5mm of n0]  {Player 2 wins};

% --- edges from 10 ---
\draw[sideedge] (n10) -- node[above left, xshift=-1mm, yshift=0.5mm] {$D(10)=10-\floor{\sqrt{10}}^2=1$} (n1);
\draw[mainedge] (n10) -- node[above right, xshift=1mm, yshift=0.5mm] {$U(10)=10+\floor{\sqrt{10}}^2=19$} (n19);

% --- edges from 19 ---
\draw[mainedge] (n19) -- node[above left, xshift=-1mm, yshift=0.5mm] {$D(19)=3$} (n3);
\draw[sideedge] (n19) -- node[above right, xshift=1mm, yshift=0.5mm] {$U(19)=35$} (n35);

% --- edges from 3 ---
\draw[sideedge] (n3) -- node[above left, xshift=-1mm, yshift=0.5mm] {$D(3)=2$} (n2);
\draw[mainedge] (n3) -- node[above right, xshift=1mm, yshift=0.5mm] {$U(3)=4$} (n4);

% --- edges from 4 ---
\draw[mainedge] (n4) -- node[above left, xshift=-1mm, yshift=0.5mm] {$D(4)=0$} (n0);
\draw[sideedge] (n4) -- node[above right, xshift=1mm, yshift=0.5mm] {$U(4)=8$} (n8);

% --- chosen sequence box ---
\node[
    draw=blue!65!black,
    rounded corners=2pt,
    line width=0.8pt,
    fill=blue!3,
    inner xsep=8pt,
    inner ysep=5pt,
    anchor=north
] at (1.15,-8.75)
{\small \emph{Chosen sequence:} $10 \to 19 \to 3 \to 4 \to 0$};

\end{tikzpicture}
\caption{A sample play of the game starting at $n=10$. The highlighted path shows the sequence
$10\to 19\to 3\to 4\to 0$, in which Player 2 wins.}
\label{fig:game-tree-10}
\credit{I specified the game and desired illustration, and Astra generated a visual design and TikZ implementation. I manually adjusted the positioning and styling of nodes and edges to improve clarity and decreased the size of the illustration. Produced in 2026.}
\end{figure}
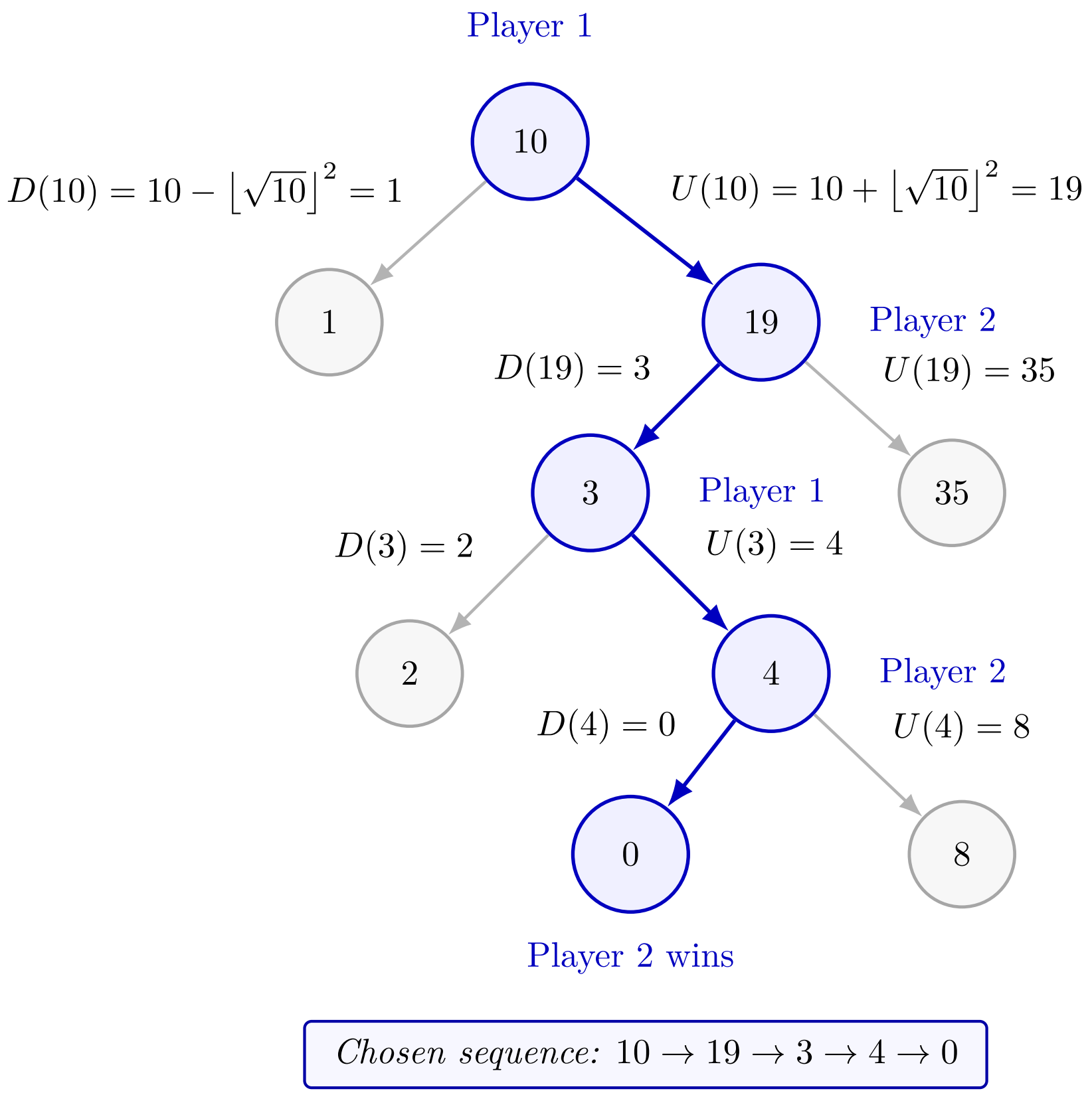

% DIAGRAM: A diagram illustrating the moves from a given position $N$ could be helpful here. For example, a tree diagram showing $N$ branching to $D(N)$ and $U(N)$, and subsequent moves from those positions.

% We now define the $\mathcal{P}$-positions and $\mathcal{N}$-positions of the game, which are typical concepts in combinatorial game theory.

% \begin{definition}
%     The \vocab{$\mathcal{P}$-positions} are losing for the player whose turn it is given optimal play.
% \end{definition}

% \begin{definition}
%     The \vocab{$\mathcal{N}$-positions} are winning for the player whose turn it is given optimal play.
% \end{definition}

% Formally,
% \begin{align*}
%     \mathcal P &= \{g \mid \forall h \in \Opt(g), h \in \mathcal N\} \\
%     \mathcal N &= \{g \mid \exists h \in \Opt(g), h \in \mathcal P\}
% \end{align*}

\begin{definition}
  A \vocab{remoteness function} is a function $R(n)$ that gives the number of
  moves until the game ends, given optimal play from position $n$.
\end{definition}
Note that for our game, $R(n) = 0$ if and only if $n = 0$, and that $R(n) = 1$ if and only if
$n$ is a positive perfect square. Some numbers, such as 2 and 3, have infinite remoteness, because two perfectly rational players would just alternate between 2 and 3 forever.

% \begin{example}[Examples of $\mathcal{P}$ and $\mathcal{N}$ positions]
%     We demonstrate by analyzing the positions $5$ and $9$.
%     \begin{itemize}
%         \item The number 5 is the smallest $\mathcal{P}$-position. The moves from 5 are $\Opt(5) = \{D(5), R(5)\} = \{1, 9\}$. Both 1 and 9 are perfect squares, so 5 is a $\mathcal{P}$-position.
%         \item The number 9 is a $\mathcal{N}$-position. The moves from 9 are $\Opt(9) = \{D(9), R(9)\} = \{0, 18\}$. Since one of these moves allow a player to win, 9 is a $\mathcal{N}$-position.
%     \end{itemize}
% \end{example}

% The first few $\mathcal{P}$-positions are $0, 5, 20, 29, 45$ \footnote{\url{https://oeis.org/A005240}}, and the first few $\mathcal{N}$-positions are $1, 4, 9, 11, 14$ \footnote{\url{https://oeis.org/A005241}}.

\begin{definition}
  An \vocab{E. P. Stein number}\footnote{Named after Richard A. Epstein's pseudonym, E. P. Stein, as per \cite{volkert}.} is a number $n$ such that $R(n) = 2$. In other
  words, an E. P. Stein number is a number such that no matter which move the
  current player takes from that number, the next turn is guaranteed to have an
  available move to $0$---i.e., a number for which there exist positive
  integers $a$ and $b$ such that $D(n) = a^2$ and $U(n) = b^2$.
\end{definition}
\begin{example}
  The number $5$ is an E. P. Stein number, since $D(5) = 1^2$ and $U(5) = 3^2$.
  Although $29$ is a losing position given optimal play for the current player,
  it is not an E. P. Stein number, since
  $U(29) = 29 + \floor{\sqrt{29}}^2 = 54$; and indeed, $R(29) = 4$.
\end{example}

The first few E. P. Stein numbers are $5, 20, 45, 80,
145$ \cite{oeisA397453}.

\begin{definition}
  The \vocab{E. P. Stein counting function} $S_2(X)$ is defined as the number of
  E. P. Stein numbers less than or equal to $X$:
  \[
    S_2(X) = \#\{n \leq X : R(n) = 2\}.
  \]
\end{definition}

The purpose of this paper is to prove the following asymptotic:

\begin{theorem}\label{thm:main}
  As $X \to \infty$, we have
  \[
    S_2(X) = \frac{1}{\pi^2}X^{1/4}\log X + O(X^{1/4}).
  \]
  Equivalently, $S_2(X) \sim \frac{1}{\pi^2}X^{1/4}\log X$; that is,
  \[
    \lim_{X \to \infty} \frac{S_2(X)}{X^{1/4}\log X} = \frac{1}{\pi^2},
  \]
  with relative error $O\left((\log X)^{-1}\right)$.
\end{theorem}

Steuding and Volkert \cite{volkert} prove that there are infinitely many E. P. Stein numbers, and develop a parameterization for them along the way. For our proof of \Cref{thm:main}, we build on their method to develop an alternate parameterization that bijectively characterizes the E. P. Stein numbers, using triples $(d, u, v)$.

After parameterizing E. P. Stein numbers in this way, we count the number of such triples satisfying the necessary inequalities using lattice geometry and equidistribution techniques. Finally, to obtain the asymptotic for $S_2(X)$, we will sum over all possible values of $d$ using Möbius inversion and a weighted harmonic sum\footnote{
  I used ChatGPT, particularly GPT-6 Astra, for
the development of mathematical proofs and ideas, assistance
with computational verification, and the creation and
refinement of TikZ figures and plots. I developed the
research question and conjecture, principal mathematical results, and reviewed
and verified any AI-suggested mathematical arguments.
AI-assisted graphics are credited individually.}.

\section{A Bijective Parameterization of E. P. Stein Numbers}

% We present for completeness parameterization found in Volkert and Steuding's work \cite{volkert}. express paramterization in terms of game moves - talk through the connectino through moves

We first focus on the parameterization of E. P. Stein numbers, which assumes the form $n = d^2(u^4 + 4v^4)$ with certain divisibility conditions on $d$, $u$, and $v$.

Let $m=\floor{\sqrt{n}}$. For an E.~P.~Stein number, $n-m^2=a^2$ and $n+m^2=b^2$. Consequently, we have
\[
b^2-a^2=2m^2,
\qquad
n^2=m^4+(ab)^2.
\]
Steuding and Volkert \cite{volkert} proved the not-necessarily-bijective parameterization
\[
n=g^2(U^4+4V^4),\qquad m=2gUV,
\]
where \(U\) is odd and \(\gcd(U,V)=1\).

% move down to example 2.3 
\begin{definition}\label{def:mustar}
  We define $\mu^{*}(n)$ to be $1$ if $p^2 \nmid n$ for all \emph{odd} primes
  $p$, and $0$ otherwise. In other words, the odd part of $n$ is square-free; the even part may have arbitrary powers of $2$.
\end{definition}

We now convert the Steuding--Volkert parameterization into a unique form with a $\mu^*$ condition rather than a coprimality condition.

\begin{theorem}[Parameterization of E. P. Stein numbers]\label{thm:param}
  Define
  \[
    \mathcal{A} = \left\{(d, u, v) \in \mathbb{Z}_{>0}^{3} : u \text{ odd},\
      \mu^*(d) = 1,\ d(u^2 - 2v^2)^2 \leq 4uv\right\}.
  \]
  Then the map
  \[
    (d, u, v) \mapsto d^2(u^4 + 4v^4)
  \]
  is a bijection from $\mathcal{A}$ to the set of E. P. Stein numbers.
\end{theorem}

\begin{proof}
  We prove that the displayed map indeed sends every triple in $\mathcal{A}$ to an
  E. P. Stein number, and that it is surjective and injective.

  % reword to proving paramterization works
  To prove that the map is well defined, take $(d, u, v) \in \mathcal{A}$ and set
  $n = d^2(u^4 + 4v^4)$ and $m = 2duv$. We then have
  \begin{align*}
    n - m^2 &= d^2(u^4 + 4v^4) - 4d^2u^2v^2 = d^2(u^2 - 2v^2)^2 = a^2, \\
    n + m^2 &= d^2(u^4 + 4v^4) + 4d^2u^2v^2 = d^2(u^2 + 2v^2)^2 = b^2.
  \end{align*}
  It now suffices to show that $m^2 < n < (m + 1)^2$. Since $m = 2duv$, we have
  $4d^2u^2v^2 < d^2(u^4 + 4v^4)$ from $0 < d^2(u^2 - 2v^2)^2$. Note that
  $u^2 - 2v^2$ cannot be zero because $u$ is odd. Moreover, because
  $d(u^2 - 2v^2)^2 \leq 4uv$, it follows that $n - m^2 \leq 2m$, and so
  $n < (m + 1)^2$. Hence $\floor{\sqrt{n}} = m$, so that $D(n) = a^2$ and
  $U(n) = b^2$ are both positive squares and $n$ is an E. P. Stein number.

  For surjectivity, let $n$ be an E. P. Stein number. The argument of Steuding
  and Volkert \cite[p.3]{volkert} yields positive integers $g, U, V$, with $U$ odd and
  $\gcd(U,V)=1$, such that
  \[
    n = g^2(U^4 + 4V^4)
  \]
  with $g(U^2 - 2V^2)^2 \leq 4UV$. Uniquely factor $g = dk^2$, where $k$ is odd
  and $\mu^*(d) = 1$: all powers of $2$ remain in $d$, while for each odd prime
  $p$, pairs of powers of $p$ go into $k^2$, leaving exponent $0$ or $1$ in $d$.
  Now set $u = kU$ and $v = kV$. We then have
  \[
    n = d^2(k^4U^4 + 4k^4V^4) = d^2(u^4 + 4v^4)
  \]
  and
  \begin{align*}
    dk^4 (U^2 - 2V^2)^2 &\leq 4k^2UV, \\
    d(u^2 - 2v^2)^2 &\leq 4uv.
  \end{align*}
  Because $k$ and $U$ are odd, $u$ is too; and $\mu^*(d) = 1$ by the construction $g = dk^2$.
  Thus $(d, u, v) \in \mathcal{A}$, as desired; in other words, every
  E. P. Stein number is in the image of the stated map.

  To prove injectivity, suppose $(d_1, u_1, v_1), (d_2, u_2, v_2) \in
  \mathcal{A}$ and $d_1^2(u_1^4 + 4v_1^4) = d_2^2(u_2^4 + 4v_2^4) = n$. Let
  $k_i = \gcd(u_i, v_i)$ for $i = 1, 2$. Then
  $n = (d_ik_i^2)^2(U_i^4 + 4V_i^4) = g_i^2(U_i^4 + 4V_i^4)$, where
  $g_i = d_ik_i^2$, $U_i = \frac{u_i}{k_i}$, and $V_i = \frac{v_i}{k_i}$;
  moreover, $\gcd(U_i, V_i) = 1$, and $k_i$ is odd because it divides the odd
  number $u_i$.

  From the parameterization argument above, $\floor{\sqrt{n}} = m_i =
  2g_iU_iV_i$, so $m := m_1 = m_2$. Let $a = \sqrt{n - m^2} =
  g_i\abs{U_i^2 - 2V_i^2}$ and $b = \sqrt{n + m^2} = g_i(U_i^2 + 2V_i^2)$.

  % label each of the math line breaks
  We now claim that $g := g_1 = g_2$. This is because
  \[\gcd(U_i^2 - 2V_i^2, U_i^2 + 2V_i^2, 2U_iV_i) = 1,\]
  which can be seen from the coprimality of $U_i$ and $V_i$ together with the oddness of $U_i$, and so
  $\gcd(a, b, m) = g_i = g$. Next, we observe that
  \[
    \left\{ U_i^2, 2V_i^2 \right\}
    = \left\{ \frac{b+a}{2g}, \frac{b-a}{2g} \right\};
  \]
  the right-hand side depends only on $n$. Of these two integers on the left-hand side of the correspondence above, $U_i^2$ is the
  odd one and $2V_i^2$ is the even one. Thus $U := U_1 = U_2$ and
  $V := V_1 = V_2$.

  % words explaning assembly of the triple from the parameters g, U, V and the uniqueness argument

  Recall from earlier that $k_i = \gcd(u_i, v_i)$. For each $i$ we have $g = g_i = d_ik_i^2$, where $k_i$ is odd and
  $\mu^*(d_i) = 1$. This factorization is unique: for every odd prime, all pairs
  of prime factors go into $k_i^2$, while the remaining exponent $0$ or $1$
  stays in $d_i$, and the entire power of $2$ stays in $d_i$. Hence $d_1 = d_2$
  and $k_1 = k_2$; it then follows that $u_1 = k_1U = k_2U = u_2$ and
  $v_1 = k_1V = k_2V = v_2$. Thus
  \[
    (d_1,u_1,v_1) = (d_2,u_2,v_2),
  \]
  so the map is injective. Having shown that the map sends $\mathcal{A}$ into the
  E. P. Stein numbers and is both injective and surjective, we conclude that it
  is a bijection, as desired.
\end{proof}

\Cref{thm:param} allows us to iterate over positive integer triples $(d, u, v)$
satisfying the stated inequality with no coprimality conditions on $u$ and $v$ (as opposed to Steuding and Volkert's original parameterization); this lets us sum over the relevant triples more directly.

% talk about the game 
\begin{example}
    Take $d = 1, u = 3, v = 2$. Then 
    \begin{align*}
        n &= 3^4 + 4 \cdot 2^4 = 81 + 64 = 145, \\
        m &= 2 \cdot 3 \cdot 2 = 12.
    \end{align*}    
    Since $D(145) = 145 - 12^2 = 1$ and $U(145) = 145 + 12^2 = 289$, the triple $(1, 3, 2)$ produces an E.~P.~Stein number, 145, as desired. In our game, both of 145's legal moves reach positive squares, and the player's opponent can reach 0 on the next turn regardless of which move is chosen.
\end{example}

% number of elements of A whose associated game position is less than $X^4$. Words explaining relationship mathcal A object to T object
\begin{definition}\label{def:Td}
  Define $T_{d}(Y)$ as the number of positive integer pairs $(u, v)$, with $u$
  odd, satisfying $d(u^2 - 2v^2)^2 \leq 4uv$ and $d^2(u^4 + 4v^4) \leq Y^4$.
\end{definition}

Throughout, $Y$ denotes the fourth root of the cutoff, so that counting E.~P.~Stein numbers up to $X$ means taking $Y = X^{1/4}$; the condition $d^2(u^4 + 4v^4) \leq Y^4$ is then exactly the requirement that the E.~P.~Stein number constructed from $(d, u, v)$ be at most $X$. We can thus think of $T_d(Y)$ as counting, for fixed $d$, the triples $(d, u, v)$ that satisfy every condition defining $\mathcal{A}$ except $\mu^*(d) = 1$, subject to the additional constraint that the associated E.~P.~Stein number is at most $X$.

As a consequence of \Cref{thm:param}, our desired count is thus
\[
  S_2(Y^4) = \sum_{d \geq 1} \mu^*(d) T_d(Y).
\]

The factor $\mu^*(d)$ appears because we have omitted the $\mu^*(d)$ condition in the definition of $T_d(Y)$.
Note, however, that $d$ satisfies
\[
  d \leq 4uv \leq 2 \sqrt{u^4 + 4v^4} \leq \frac{2Y^2}{d},
\]
and thus $d \leq \sqrt{2}Y$; otherwise $T_d(Y) = 0$. This allows us to evaluate the sum
\[
S_2(Y^4)=\sum_{1\leq d\leq\sqrt{2}Y}\mu^*(d)T_d(Y).
\]

\section{Lattice Geometry and Equidistribution}

We have now converted the problem of counting E.~P.~Stein numbers into a number-theoretic problem of counting integer points in a region defined by certain inequalities. Moving forward, we will estimate the distribution of the lattice points in this valid region by utilizing a change of variables and approximating each lattice point count with the corresponding area of the region it occupies.

Our eventual counting problem will reduce to determining how often
the fractional parts of $v/\sqrt{2}$ lie in a short interval; the following lemma provides a rough estimate for this count.
% fix wording here about circle
\begin{lemma}\label{lem:equidistribution}
  Let $\alpha = \frac{1}{\sqrt{2}}$, and let $\frac{p}{q}$ be a rational number
  with $p, q > 0$, $\gcd(p, q) = 1$, and
  \[
    \abs{\alpha - \frac{p}{q}} < \frac{1}{q^2}.
  \]
  Then, uniformly over real phases $\theta$, over intervals $I \subseteq [0, 1)$
  of length $w \in [0, 1)$ taken modulo $1$, and over integers $V \geq 1$, we have
  \[
    n(V) = \#\{1 \leq v \leq V : \fracpart{\alpha v + \theta} \in I\}
    = Vw + O\left(\frac{V}{q} + qw + 1\right),
  \]
  where $\fracpart{x}$ denotes the fractional part of $x$, so that
  $\fracpart{\alpha v + \theta}$ is $\alpha v + \theta \pmod{1}$.
\end{lemma}

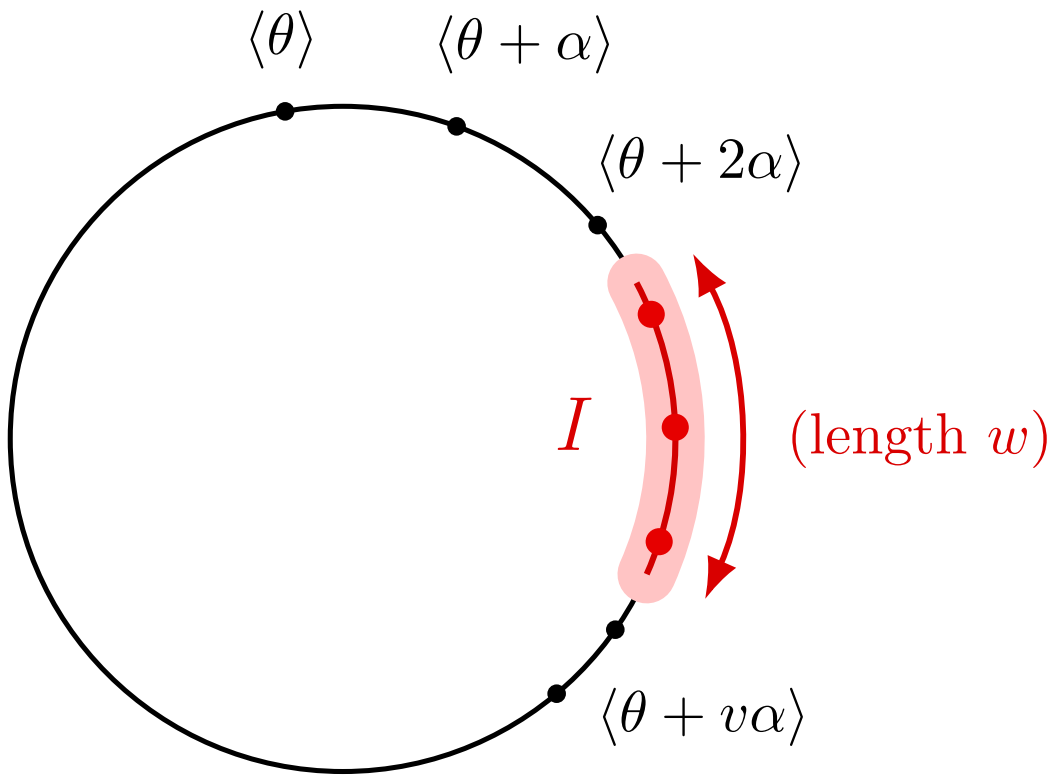
\begin{figure}[ht]
\centering
\begin{tikzpicture}[font=\normalsize, >=Latex]

\def\R{2.20}

% Unit circle
\draw[line width=0.95pt] (0,0) circle[radius=\R];

% Highlighted interval I
\draw[red!23, line width=11pt, line cap=round]
  (-24:\R) arc[start angle=-24, end angle=28, radius=\R];
\draw[red!82!black, line width=1.15pt]
  (-24:\R) arc[start angle=-24, end angle=28, radius=\R];

% More evenly spaced black orbit points
\foreach \ang in {100,310,325,40,70}{
  \fill[black] (\ang:\R) circle[radius=1.75pt];
}

% Red points lying in I
\foreach \ang in {-18,2,22}{
  \fill[red!90!black] (\ang:\R) circle[radius=2.55pt];
}

% Labels
\node[anchor=south] at (100:2.35) {$\fracpart{\theta}$};
\node[anchor=west]  at (80:2.65) {$\fracpart{\theta+\alpha}$};
\node[anchor=west]  at (50:2.38) {$\fracpart{\theta+2\alpha}$};
\node[anchor=west]  at (310:2.40) {$\fracpart{\theta+v\alpha}$};

% Curved double-headed arrow indicating arc length
\draw[red!85!black, <->, line width=1.1pt]
  (-24:2.62) arc[start angle=-24, end angle=28, radius=2.62];

% Separate labels for I and length w
\node[red!85!black, font=\Large] at (1.53,0.1) {$I$};
\node[red!85!black, anchor=west] at (2.80,0) {$(\text{length } w)$};

\end{tikzpicture}
\caption{Schematic illustration of the orbit points $\fracpart{\theta+k\alpha}$ on the unit interval---illustrated as a circle here for visualizing the modular arithmetic---and the highlighted interval $I$ of length $w$, where $n(V)$ counts the number of points falling into $I$. The plotted positions of the orbit points are schematic and spaced for legibility and are not the actual positions of $\fracpart{\theta+k\alpha}$ for any particular $\theta$ and $\alpha$.}
\label{fig:circle-orbit}
\credit{The initial diagram was generated with Astra; Astra subsequently supplied TikZ code and revisions based on my visual feedback, and then I manually edited the positioning and size of various nodes to improve readability and spacing. Produced in 2026.}
\end{figure}

%% Keep Figure 2 with the statement of the lemma, not inside its proof.
\FloatBarrier

\begin{proof}
  We define a \textit{block} to be $q$ consecutive terms of the
sequence \[\left(\fracpart{\theta+v\alpha}\right)_{v\geq 1}.\] Say that we are analyzing the block $\fracpart{\theta + (\ell + j)\alpha}$ for
  $0 \leq j < q$. Comparing this with
  $\fracpart{\theta + \ell\alpha + j\frac{p}{q}}$, continued fraction theory for non-zero $j$ \cite{cont_fractions}
  yields the inequality
  \[
    j\abs{\alpha - \frac{p}{q}} \leq \frac{j}{q^2} < \frac{1}{q},
  \]
  and so each actual point is within $\frac{1}{q}$ of the corresponding rational
  point. Since $\gcd(p,q)=1$, the rational comparison points
form a translate modulo $1$ of the evenly spaced grid
\[
  \left\{0,\frac1q,\ldots,\frac{q-1}{q}\right\}.
\]
Hence an interval of length $w$ contains $qw+O(1)$
of these comparison points. Moving a comparison point by less than $1/q$ can change
whether it belongs to $I$ only if it lies within distance
$1/q$ of an endpoint of $I$. Since the grid points are
spaced $1/q$ apart, only $O(1)$ points can lie in these
endpoint neighborhoods. Therefore, each complete block
contains $qw+O(1)$ actual orbit points in $I$, uniformly
in the starting index, phase, and interval.

Now write $V=Mq+H$, where $0\leq H<q$. The first $Mq$
terms consist of $M$ complete blocks, each contributing
$qw+O(1)$ points in $I$. Let $R_H$ denote the number of points in $I$ among the
remaining $H$ terms. Since these terms are contained in
one additional complete block, we have
\[
  0\leq R_H\leq qw+O(1).
\]
Moreover, $0\leq Hw\leq qw$, so
\[
  R_H=Hw+O(qw+1).
\]
Combining the complete blocks and the remainder gives
\begin{align*}
  n(V)
  &=M(qw+O(1))+R_H\\
  &=Mqw+Hw+O(M+qw+1)\\
  &=Vw+O\left(\frac{V}{q}+qw+1\right),
\end{align*}
as desired.
\end{proof}

We now employ a ``sandwich'' argument and a change of variables---where
$t = u - \sqrt{2}v$---to trap the admissible pairs $(t, v)$ between two
rectangles: an inner rectangle $\mathcal{R}_{\mathrm{in}}$ and an outer
rectangle $\mathcal{R}_{\mathrm{out}}$. Define
\[
  \Omega_d(Y) = \left\{(t,v) :
  \begin{array}{l}
    v \ge 1,\ u = t + \sqrt2 v > 0,\\
    d(u^2-2v^2)^2 \le 4uv,\\
    d^2(u^4+4v^4) \le Y^4
  \end{array}
  \right\}.
\]

% swap v and t back; admissable diagram is shitty wording
Note that $t = u - \sqrt{2}v$ is in general irrational, so the points of $\Omega_d(Y)$ are not themselves lattice points. Rather, $T_d(Y)$ counts the integer pairs $(u,v)$, with $u$ odd, whose transformed coordinates $(t,v)$ lie in $\Omega_d(Y)$.
The difference between the two rectangles is that if
$(t,v) \in \mathcal{R}_{\mathrm{in}}$, it necessarily follows that
$(t,v) \in \Omega_{d}(Y)$; whereas if $(t,v) \in \Omega_d(Y)$, it only follows
that $(t,v) \in \mathcal{R}_{\mathrm{out}}$. The next three propositions
establish precise bounds for $(t,v)$ within these rectangles, and hence bound
the $\Omega_d(Y)$ region; in other words, they show that
$\mathcal{R}_{\mathrm{in}} \subseteq \Omega_{d}(Y) \subseteq
\mathcal{R}_{\mathrm{out}}$.

% make caption a little clearer and show sizing; can do it in area terms
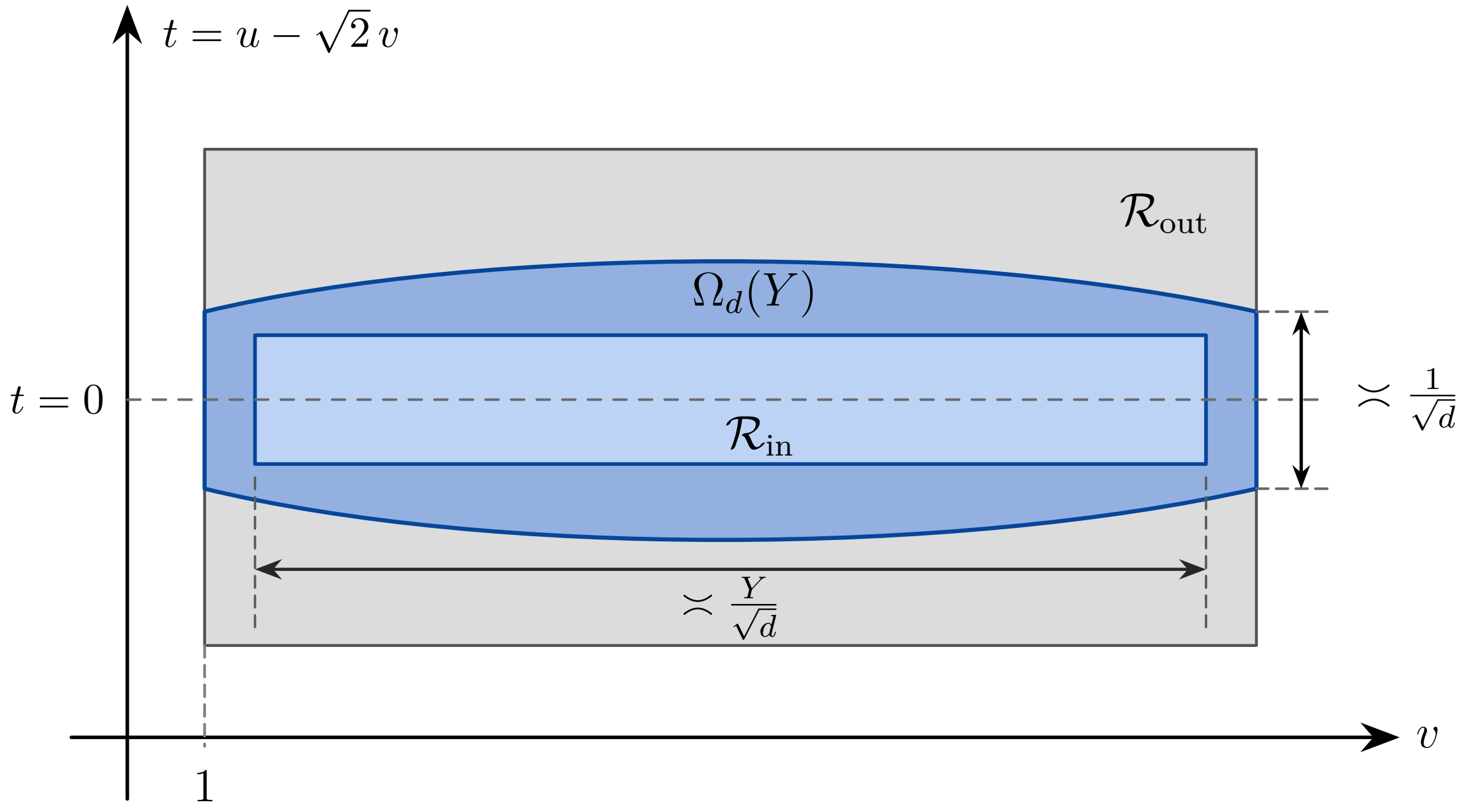
\begin{figure}[ht]
    \centering
    \resizebox{0.98\linewidth}{!}{%
        % Requires \usepackage{tikz} and \usetikzlibrary{arrows.meta,calc,shadings}.
\definecolor{outergray}{RGB}{220,220,220}
\definecolor{bandblue}{RGB}{147,176,224}
\definecolor{innerblue}{RGB}{189,211,246}
\definecolor{strokeBlue}{RGB}{6,70,155}

\begin{tikzpicture}[x=1cm,y=1cm,
  >={Stealth[length=3.0mm,width=2.0mm]},
  line cap=round,line join=round]

% Title and axis label (separate zones to avoid collisions)
\node[anchor=west,font=\normalsize] at (0.18,6.31) {$t=u-\sqrt{2}\,v$};

% Coordinate axes: the horizontal axis is below the displaced strip
\draw[-{Stealth[length=3.5mm]},line width=0.90pt]
  (-0.50,0) -- (11.36,0) node[right,font=\large] {$v$};
\draw[-{Stealth[length=3.5mm]},line width=0.90pt]
  (0,-0.55) -- (0,6.54);

% Outer comparison rectangle: flat gray field
\fill[outergray]
  (0.69,0.82) rectangle (10.08,5.25);
\draw[black!68,line width=0.72pt]
  (0.69,0.82) rectangle (10.08,5.25);
\node[anchor=north east,font=\large] at (9.80,4.99) {$\mathcal{R}_{\mathrm{out}}$};

% Admissible region: a very shallow, smoothly bulging thin strip.
% Schematic: shape is NOT an equation-by-equation plot of the boundary.
\fill[bandblue]
  (0.69,3.80)
  .. controls (3.20,4.40) and (7.37,4.40) .. (10.08,3.80)
  -- (10.08,2.22)
  .. controls (7.37,1.61) and (3.20,1.61) .. (0.69,2.22)
  -- cycle;
\draw[strokeBlue,line width=1.05pt]
  (0.69,3.80)
  .. controls (3.20,4.40) and (7.37,4.40) .. (10.08,3.80)
  -- (10.08,2.22)
  .. controls (7.37,1.61) and (3.20,1.61) .. (0.69,2.22)
  -- cycle;

% Inner comparison rectangle, clearly entirely inside the blue region
\fill[innerblue] (1.14,2.44) rectangle (9.63,3.59);
\draw[strokeBlue,line width=0.93pt]
  (1.14,2.44) rectangle (9.63,3.59);
\node[font=\large] at (5.60,3.96) {$\Omega_d(Y)$};
\node[font=\large] at (5.65,2.69) {$\mathcal{R}_{\mathrm{in}}$};

% Centerline (on top of all shaded fills), at t=0
\draw[black!58,dash pattern=on 3.3pt off 3.2pt,line width=0.65pt]
  (0.00,3.015) -- (10.65,3.015);
\node[anchor=east,font=\normalsize] at (-0.07,3.015) {$t=0$};

% Left edge of the outer domain corresponds to v=1
\draw[black!50,densely dashed,line width=0.7pt] (0.69,0.82) -- (0.69,-0.08);
\node[font=\large,anchor=north] at (0.69,-0.17) {$1$};

% Horizontal scale of the inner comparison rectangle
\draw[black!65,densely dashed,line width=0.60pt]
  (1.14,2.32) -- (1.14,0.98);
\draw[black!65,densely dashed,line width=0.60pt]
  (9.63,2.32) -- (9.63,0.98);
\draw[{Stealth[length=2.4mm]}-{Stealth[length=2.4mm]},black!85,line width=0.80pt]
  (1.14,1.5) -- (9.63,1.5);
\node[font=\normalsize,inner sep=0.6pt] at (5.39,1.15)
  {$\asymp\frac{Y}{\sqrt d}$};

% Vertical scale bracket with adequate right margin
\draw[black!60,densely dashed,line width=0.65pt]
  (10.08,3.80) -- (10.72,3.80);
\draw[black!60,densely dashed,line width=0.65pt]
  (10.08,2.22) -- (10.72,2.22);
\draw[{Stealth[length=2.1mm]}-{Stealth[length=2.1mm]},black,line width=0.85pt]
  (10.48,2.22) -- (10.48,3.80);
\node[anchor=west,font=\large] at (10.83,3.015)
  {$\asymp\frac{1}{\sqrt d}$};

\end{tikzpicture}
    }
    \caption{Schematic illustration of
    $\mathcal{R}_{\mathrm{in}}\subseteq\Omega_d(Y)
    \subseteq \mathcal{R}_{\mathrm{out}}$. We use this schematic to visualize bounding the area of $\Omega_d(Y)$ between the areas of the inner and outer rectangles.}
    \label{fig:thin-strip}
    \credit{Astra generated the original visual, then recreated it in TikZ after I requested the design. The figure is schematic rather than an exact plot of \(\Omega_d(Y)\)---I manually adjusted the shape of some of the boundaries and the positioning of the rectangles and arrows for readability. Produced in 2026.}
\end{figure}

Throughout, we use the conventions $c = 2^{-1/4}$, $c' = 2^{-3/4}$, and
$Y^4 = X$.

The inequality $d(u^2-2v^2)^2\le 4uv$ defining $\Omega_d(Y)$ is nonlinear in $u$
and $v$, but becomes simpler when expressed in terms of the
distance $t=u-\sqrt{2}v$ from the line $u=\sqrt{2}v$.
The following equivalence allows us to bound the width of
$\Omega_d(Y)$ in the $t$-direction.

%% Keep Figure 3 ahead of the propositions it motivates, so it cannot be
%% deferred into the middle of one of their proofs.
\FloatBarrier

% for each proposition, give some motivation
\begin{proposition}\label{prop:criterion}
  Let $d \ge 1$, $u > 0$, and $v \ge 1$, and let $t = u - \sqrt2 v$ and
  $s = u + \sqrt2 v$. Then
  \[
    d(u^2-2v^2)^2 \le 4uv
    \quad\Longleftrightarrow\quad
    t^2\left(d+\frac1{\sqrt2 s^2}\right) \le \frac1{\sqrt2}.
  \]
\end{proposition}

\begin{proof}
  First, note that $s$ and $t$ satisfy $s^2 t^2 = (u^2 - 2v^2)^2$ and
  $\frac{s^2 - t^2}{\sqrt{2}} = 4uv$. Using the condition
  $d(u^2 - 2v^2)^2 \leq 4uv$, we find that
  \begin{align*}
    d s^2 t^2 &\leq \frac{s^2 - t^2}{\sqrt{2}} \\
    d t^2 &\leq \frac{1}{\sqrt{2}} - \frac{t^2}{s^2\sqrt{2}} \\
    t^2 \left(d + \frac{1}{\sqrt{2}s^2}\right) &\leq \frac{1}{\sqrt{2}},
  \end{align*}
  as desired. The converse is obtained by reversing these steps.
\end{proof}

To obtain a lower bound for $T_d(Y)$, we first construct a
rectangle contained entirely in $\Omega_d(Y)$.
Counting integer pairs in this rectangle will provide a lower bound for the count in $\Omega_d(Y)$.

\begin{proposition}\label{prop:inner}
  For $d \geq 1$ and $Y \geq 2$, let
  \[
    \mathcal{R}_{\mathrm{in}} = \left\{(t,v) : 1 \leq v \leq
      \frac{c'(Y - c)}{\sqrt{d}},\
      \abs{t} \leq \frac{c}{\sqrt{d + \frac{1}{2\sqrt{2}}}}\right\}.
  \]
  Then $\mathcal{R}_{\mathrm{in}} \subseteq \Omega_{d}(Y)$.
\end{proposition}

\begin{proof}
  Suppose $(t,v) \in \mathcal{R}_{\mathrm{in}}$.
  Apply the $L^4$ norm to the points $(u, \sqrt{2}v)$ and
  $(\sqrt{2}v, \sqrt{2}v)$, together with the triangle inequality, to obtain
  \begin{align*}
    \norm{(u, \sqrt{2}v)}_4 &\leq \norm{(u - \sqrt{2}v, 0)}_4
      + \norm{(\sqrt{2}v, \sqrt{2}v)}_4 \\
    (u^4 + 4v^4)^{1/4} &\leq \abs{t} + 2^{3/4} v,
  \end{align*}
  where $t = u - \sqrt{2}v$ as before. Since $v$ satisfies the bound defining
  $\mathcal{R}_{\mathrm{in}}$,
  \[
    (u^4 + 4v^4)^{1/4} \leq \abs{t} + 2^{3/4} v
    < \frac{Y - c}{\sqrt{d}} + \frac{c}{\sqrt{d}} = \frac{Y}{\sqrt{d}}.
  \]
  Rearranging shows that $d^2(u^4 + 4v^4) \leq Y^4 = X$, as desired. Next, since
  $v \geq 1$ and $\abs{t} \leq c/\sqrt{d + \frac{1}{2\sqrt{2}}} < \sqrt{2}$, we
  have $u = t + \sqrt{2}v > 0$; hence $s = u + \sqrt{2}v > \sqrt{2}v \geq \sqrt{2}$,
  so $s^2 \geq 2$ and
  \[
    t^2 \left(d + \frac{1}{\sqrt{2}s^2}\right)
    \leq t^2\left(d + \frac{1}{2\sqrt{2}}\right) \leq \frac{1}{\sqrt{2}},
  \]
  which by \Cref{prop:criterion} is the required inequality.
\end{proof}

For the corresponding upper bound, we enclose $\Omega_d(Y)$
in a larger rectangle. Together with \Cref{prop:inner}, this will allow us to estimate $T_d(Y)$ by counting points in
two simpler regions.

\begin{proposition}\label{prop:outer}
  For $d \geq 1$ and $Y \geq 2$, let
  \[
    \mathcal{R}_{\mathrm{out}} = \left\{(t,v) : 1 \leq v \leq
      \frac{c'(Y + c)}{\sqrt{d}},\ \abs{t} \leq \frac{c}{\sqrt{d}}\right\}.
  \]
  Then $\Omega_{d}(Y) \subseteq \mathcal{R}_{\mathrm{out}}$.
\end{proposition}

\begin{proof}
  Suppose that $(t,v) \in \Omega_{d}(Y)$. By \Cref{prop:criterion}, we know
  that $t^2\left(d + \frac{1}{\sqrt{2}s^2}\right) \leq \frac{1}{\sqrt{2}}$, and
  thus $t^2 d \leq \frac{1}{\sqrt{2}}$. Hence
  $\abs{t} \leq \frac{1}{\sqrt[4]{2}\sqrt{d}} = \frac{c}{\sqrt{d}}$, as desired.

  % same l4 argument to prove inequality in proposition 3.3
  From the same $L^4$ norm argument we used to prove the inequality in \Cref{prop:inner}, we have
  \[
    2^{3/4}v \leq \abs{t} + (u^4 + 4v^4)^{1/4},
  \]
  and thus
  \[
    2^{3/4}v \leq \frac{Y}{\sqrt{d}} + \frac{c}{\sqrt{d}},
  \]
  so that $v \leq \frac{c'(Y + c)}{\sqrt{d}}$, as desired.
\end{proof}

Now that we have inequalities describing the inner and outer rectangles that sandwich $\Omega_d(Y)$, we can estimate the number of integer pairs $(u,v)$ with odd $u$ in these regions to bound $T_d(Y)$.

% transition sentence
\begin{lemma}\label{lem:Td}
  The quantity $T_d(Y)$ satisfies
  \[
    T_d(Y) = \frac{Y}{2d}
    + O\left(\frac{Y}{d^2} + \sqrt{\frac{Y}{d}} + 1\right).
  \]
\end{lemma}

% change q in previous passage
\begin{proof}
  Let $N_{\pm}$ count the number of integer pairs with odd $u$ in our
  rectangles, with $N_{+}$ counting those in the outer rectangle
  $\mathcal{R}_{\mathrm{out}}$ and $N_{-}$ counting those in the inner rectangle
  $\mathcal{R}_{\mathrm{in}}$. Notice that counting the odd $u$ satisfying
  $\abs{u - \sqrt{2}v} \leq h < 1$ is equivalent to counting the solutions of
  \begin{align*}
    \abs{2k + 1 - \sqrt{2}v} &\leq h \\
    \abs{k + \frac{1}{2} - \frac{1}{\sqrt{2}}v} &\leq \frac{h}{2} \\
    -\frac{h}{2} \leq \frac{v}{\sqrt{2}} - k - \frac{1}{2} &\leq \frac{h}{2}.
  \end{align*}
  Thus the fractional part of $\frac{v}{\sqrt{2}}$ must lie in the interval
  \[
    \fracpart{\frac{v}{\sqrt{2}}}
    \in \left[\frac{1}{2} - \frac{h}{2}, \frac{1}{2} + \frac{h}{2}\right].
  \]
  Both the inner and outer rectangles have the form $1 \le v \le L$ and
  $\abs{u - \sqrt2 v} \le h$, where $u$ must be odd. Thus we may apply
  \Cref{lem:equidistribution} to obtain
  \[
    N_{\pm} = \floor{L_{\pm}} h_{\pm}
    + O \left(\frac{\floor{L_{\pm}}}{q} + qh_{\pm} + 1\right).
  \]

  % standard recursion for Pell's equations
  Now select a convergent denominator $q$ such that
  $\sqrt{Y} \leq q \leq 6\sqrt{Y}$, which is possible since the denominators
  increase by a factor less than $6$ by a standard recursion in Pell's equations \cite{cont_fractions}; this gives the
  approximation $q \asymp \sqrt Y$. From \Cref{prop:inner,prop:outer}, we know
  that $L_{\pm} = O\left(\frac{Y}{\sqrt{d}}\right)$ and
  $h_{\pm} = O\left(\frac{1}{\sqrt{d}}\right)$.

  We now evaluate the error terms in the expression for $N_{\pm}$. Since
  $q \geq \sqrt{Y}$, we have
  $\frac{\floor{L_{\pm}}}{q} \leq \frac{L_{\pm}}{q}
  = O\left(\frac{Y/\sqrt{d}}{\sqrt{Y}}\right)
  = O\left(\sqrt{\frac{Y}{d}}\right)$ and
  $qh_{\pm} \leq 6 \sqrt{Y} \frac{c}{\sqrt{d}}
  = O\left(\sqrt{\frac{Y}{d}}\right)$. Thus the overall error term simplifies to
  \[
    N_{\pm} = L_{\pm} h_{\pm} + O\left(\sqrt{\frac{Y}{d}} + 1\right).
  \]
  We then have
  \begin{align*}
    L_{+} h_{+} &= \frac{c'(Y + c)}{\sqrt{d}} \cdot \frac{c}{\sqrt{d}}
      = \frac{c c'(Y + c)}{d}
      = \frac{Y + c}{2d} = \frac{Y}{2d} + O(1),
  \end{align*}
  and similarly
  \begin{align*}
    L_{-} h_{-} &= \frac{c'(Y - c)}{\sqrt{d}}
      \cdot \frac{c}{\sqrt{d + \frac{1}{2\sqrt{2}}}}
      = \frac{Y - c}{2\sqrt{d\left(d + \frac{1}{2\sqrt{2}}\right)}} \\
    &= \frac{Y}{2d} + O\left(\frac{Y}{d^2} + 1\right).
  \end{align*}
  Writing $E_d(Y) = \frac{Y}{d^2} + \sqrt{\frac{Y}{d}} + 1$, the preceding
  estimates give $N_{\pm} = \frac{Y}{2d} + O(E_d(Y))$.
  Since $\mathcal{R}_{\mathrm{in}} \subseteq \Omega_d(Y) \subseteq
  \mathcal{R}_{\mathrm{out}}$, the counts satisfy $N_{-} \leq T_d(Y) \leq N_{+}$,
  and we obtain
  \[
    T_d(Y) = \frac{Y}{2d}
    + O\left(\frac{Y}{d^2} + \sqrt{\frac{Y}{d}} + 1\right),
  \]
  as desired.
\end{proof}

\section{Proof of the Asymptotic}

Having developed our machinery, we now turn to the proof of \Cref{thm:main},
which we restate here for convenience: as $X \to \infty$, the counting function
of the E. P. Stein numbers---that is, of the losing positions of remoteness
$2$---satisfies
\[
  S_2(X) = \frac{1}{\pi^2}X^{1/4}\log X + O(X^{1/4}),
\]
so that the relative error is $O\left((\log X)^{-1}\right)$.

% arises from dirichlet sum that comes from paramterization

The constant $\frac{1}{\pi^2}$ is of particular interest: as we will see, it
arises from a Dirichlet sum that comes from the parameterization of every E. P. Stein number.

In \Cref{lem:mobiussum} below it will be convenient to write $K = \sqrt{2} Y$.
% explain a little more

\begin{lemma}\label{lem:mobiussum}
  We have
  \[
    \sum_{d \leq K} \frac{\mu^*(d)}{d} = \frac{8}{\pi^2} \log K + O(1).
  \]
\end{lemma}

The proof rests on the following elementary identity.

\begin{proposition}\label{prop:mustar}
  The odd square-free indicator function $\mu^{*}(d)$ is given by
  \[
    \mu^{*}(d) = \sum_{k^2 \mid d, \, k \text{ odd}} \mu (k),
  \]
  where $\mu$ is the classical Möbius function.
\end{proposition}

\begin{proof}
  The sum is equivalent to
  \[
    \mu^{*}(d) = \prod_{p^2 \mid d, \, p \text{ odd}} (1 + \mu(p)).
  \]
  This is exactly the definition of the odd square-free indicator function
  $\mu^{*}(d)$, since if any one of the odd primes squared divides $d$ then the
  corresponding term in the product is 0, and otherwise the product is 1.
\end{proof}

\begin{proof}[Proof of \Cref{lem:mobiussum}]
  We apply \Cref{prop:mustar} to evaluate the sum:
  \begin{align*}
    \sum_{d \leq K} \frac{\mu^*(d)}{d}
      &= \sum_{d \leq K} \frac{1}{d} \sum_{k^2 \mid d, \, k \text{ odd}} \mu (k) \\
      &= \sum_{k \text{ odd}} \mu(k) \sum_{\substack{d \leq K \\ k^2 \mid d}} \frac{1}{d} \\
      &= \sum_{\substack{k\le\sqrt K\\k\text{ odd}}} \frac{\mu(k)}{k^2} \sum_{m \leq K/k^2} \frac{1}{m} \\
      &= \sum_{\substack{k\le\sqrt K\\k\text{ odd}}} \frac{\mu(k)}{k^2} \left( \log \frac{K}{k^2} + O(1) \right) \\
      &= \log K \sum_{\substack{k\le\sqrt K\\k\text{ odd}}} \frac{\mu(k)}{k^2}
        - 2 \sum_{\substack{k\le\sqrt K\\k\text{ odd}}} \frac{\mu(k) \log k}{k^2} + O(1) \\
      &= \frac{4}{3} \cdot \frac{6}{\pi^2} \log K + O(1)
        = \frac{8}{\pi^2} \log K + O(1). \qedhere
  \end{align*}
\end{proof}

\begin{lemma}\label{lem:errorsum}
  We have
  \[
    \sum_{d \leq \sqrt{2}Y}
    \left(\frac{Y}{d^2} + \sqrt{\frac{Y}{d}} + 1\right) = O(Y).
  \]
\end{lemma}
\begin{proof}
  We split the sum into three separate sums:
  \begin{align*}
    \sum_{d \leq \sqrt{2}Y} \left(\frac{Y}{d^2} + \sqrt{\frac{Y}{d}} + 1\right)
      &= \sum_{d \leq \sqrt{2}Y} \frac{Y}{d^2}
        + \sum_{d \leq \sqrt{2}Y} \sqrt{\frac{Y}{d}}
        + \sum_{d \leq \sqrt{2}Y} 1 \\
      &\leq Y \sum_{d \leq \sqrt{2}Y} \frac{1}{d^2}
        + \sqrt{Y} \sum_{d \leq \sqrt{2}Y} \frac{1}{\sqrt{d}}
        + \sum_{d \leq \sqrt{2}Y} 1 \\
      &\leq Y \cdot O(1) + \sqrt{Y} \cdot O(\sqrt{Y}) + O(Y) \\
      &= O(Y) + O(Y) + O(Y) = O(Y). \qedhere
  \end{align*}
\end{proof}

We are now in a position to prove the main theorem.

\begin{proof}[Proof of \Cref{thm:main}]
  Substituting $K = \sqrt{2} Y$ in \Cref{lem:mobiussum} and combining it with
  \Cref{lem:Td,lem:errorsum} yields
  \[
    \begin{aligned}
      S_2(Y^4)
      &= \sum_{d\le\sqrt2Y}\mu^*(d)T_d(Y)\\
      &= \frac Y2\sum_{d\le\sqrt2Y}\frac{\mu^*(d)}d
        + O\left(\sum_{d\le\sqrt2Y}
          \left[\frac Y{d^2}+\sqrt{\frac Yd}+1\right]\right)\\
      &= \frac Y2\left(\frac8{\pi^2}\log(\sqrt2Y)+O(1)\right)+O(Y)\\
      &= \frac4{\pi^2}Y\log Y+O(Y) \\
      &= \frac{1}{\pi^2} X^{1/4} \log X + O(X^{1/4}),
    \end{aligned}
  \]
  as desired. Here we used that $0 \le \mu^*(d) \le 1$, so the weighted errors
  are at most the unweighted ones.
\end{proof}

\section{Computational Results}

% empirical computation up to $X = 10^{28}$
While \Cref{thm:main} gives an asymptotic for $S_2(X)$, it says nothing about how
closely that asymptotic tracks $S_2(X)$ at any particular $X$. Thus here we perform an
exact enumeration of the E. P. Stein numbers up to $X = 10^{28}$, so that the
asymptotic can be checked against exact counts over a wide range. In our proof earlier, we counted normalized triples $(d, u, v)$ with $\mu^*(d) = 1$ and no coprimality restrictions on $u$ and $v$. For the enumeration below, however, we instead iterate over coprime pairs $(U, V)$ with $U$ odd, where $g$ denotes the scaling factor in the original Steuding--Volkert parameterization.

For any triple $(g, U, V)$ and an E. P. Stein number $n$, we have
\begin{align*}
  n = g^2 (U^4 + 4V^4) &\leq X \\
  g(U^2 - 2V^2)^2 &\leq 4UV.
\end{align*}
Given a coprime pair $(U, V)$, then, we can count the number of valid $g$ to sum over as

\[
G_X(U,V)=
\min\left\{
\left\lfloor\sqrt{\frac{X}{U^4 + 4V^4}}\right\rfloor,
\left\lfloor\frac{4UV}{(U^2 - 2V^2)^2}\right\rfloor
\right\}.
\]

To sum over $U$ and $V$, notice that $4V^4 \leq U^4 + 4V^4 \leq X$, and so it suffices to consider

\[
1\le V\le
\left\lfloor\left(\frac X4\right)^{1/4}\right\rfloor.
\]

By \Cref{prop:criterion}, any admissible $U$ satisfies $\abs{U - \sqrt{2}V} \leq c = 2^{-1/4}$, so it lies in an interval of length $2c < 2$; since consecutive odd numbers differ by $2$, at most one odd \(U\) can satisfy this condition for a given \(V\). Our program thus iterates over all valid values of $V$ and, for each $V$, checks the unique odd $U$ (if it exists) that satisfies the conditions. It then sums $G_X(U,V)$ over all such pairs to compute the exact count of E. P. Stein numbers. At \(X=10^{28}\), the bound is only 7,071,067 candidate values of \(V\), making the computation feasible; the results are summarized in \Cref{tab:stein-asymptotic}.

\begin{table}[ht]
\centering
\caption{Exact counts of E. P. Stein numbers compared with the leading
asymptotic $\pi^{-2}X^{1/4}\log X$.}
\label{tab:stein-asymptotic}
\begin{tabular}{r r r r}
\toprule
$X$ & Exact $S_2(X)$ & $\pi^{-2}X^{1/4}\log X$ & Ratio \\
\midrule
$10^{4}$ & 15 & 9.3 & 1.607 \\
$10^{8}$ & 281 & 186.6 & 1.506 \\
$10^{12}$ & 3,329 & 2,799.6 & 1.189 \\
$10^{16}$ & 42,367 & 37,328.1 & 1.135 \\
$10^{20}$ & 543,144 & 466,601.3 & 1.164 \\
$10^{24}$ & 6,178,005 & 5,599,215.5 & 1.103 \\
$10^{28}$ & 70,201,703 & 65,324,181.2 & 1.075 \\
\bottomrule
\end{tabular}
\credit{Table 1 created by the author, Cameron Rampell; the exact counts were
computed with the author's own Python implementation of the parameterization
in Section 2. Produced in 2026.}
\end{table}

When comparing the ratio of the exact counts to the predicted asymptotic values, we observe that the ratio tends to 1 as \(X\) increases, but it does so in a nonmonotonic, oscillatory fashion. This suggests that by capturing the oscillations, an error term of order $(\log X)^{-1}$ may be best possible for the asymptotic formula.

\pgfplotstableread{
k exact predicted ratio
4 15 9.33202588 1.60736802
5 36 20.74368685 1.73546777
6 61 44.26568546 1.37804259
7 125 91.83621652 1.36111879
8 281 186.64051765 1.50556805
9 470 373.38636331 1.25874977
10 903 737.76142430 1.22397291
11 1939 1443.14054539 1.34359748
12 3329 2799.60776471 1.18909514
13 6301 5393.35858119 1.16828872
14 12995 10328.65994018 1.25814966
15 22706 19679.18925527 1.15380770
16 42367 37328.10352950 1.13498935
17 84784 70528.53529245 1.20212336
18 150238 132797.05637376 1.13133532
19 277478 249269.73056674 1.11316364
20 543144 466601.29411872 1.16404306
21 971449 871234.84773027 1.11502542
22 1781430 1623075.13345701 1.09756472
23 3429334 3017475.68580790 1.13649101
24 6178005 5599215.52942462 1.10336974
25 11247974 10371843.42536040 1.08447202
26 21403312 19181797.03176462 1.11581370
28 70201703 65324181.17662056 1.07466641
}\steinData

\pgfplotsset{
  epsteinaxis/.style={
    width=0.94\linewidth,
    height=7.0cm,
    xmin=4, xmax=28,
    xtick={4,8,12,16,20,24,28},
    xlabel={$k=\log_{10} X$},
    tick align=outside,
    tick style={black!50},
    grid=major,
    grid style={black!12},
    line width=0.7pt,
    legend cell align=left,
    legend style={draw=none,fill=none,font=\footnotesize},
    label style={font=\small},
    tick label style={font=\footnotesize},
  }
}

% % GRAPH 1: Exact counts and leading asymptotic, logarithmic vertical axis.
% \begin{figure}[ht!]
% \centering
% \begin{tikzpicture}
%   \begin{axis}[
%     epsteinaxis,
%     ymode=log,
%     log basis y=10,
%     ymin=8, ymax=100000000,
%     ylabel={$S_2(10^k)$},
%     legend pos=north west,
%   ]
%     \addplot[black,very thick,mark=*,mark size=1.15pt]
%       table[x=k,y=exact]{\steinData};
%     \addlegendentry{Exact count $S_2(X)$}

%     \addplot[black!60,thick,densely dashed,mark=none]
%       table[x=k,y=predicted]{\steinData};
%     \addlegendentry{Asymptotic $X^{1/4}\log X/\pi^2$}
%   \end{axis}
% \end{tikzpicture}
% \caption{Exact counts of E. P. Stein numbers compared with the leading
% asymptotic, evaluated at $X=10^k$. The vertical axis is logarithmic.}
% \label{fig:epstein-exact-asymptotic}
% \credit{Graph created by the author, Cameron Rampell, from the counts in
% Table 1; plotted with \textsf{pgfplots}.}
% \end{figure}

% GRAPH 2: Ratio to asymptotic (should tend to 1).
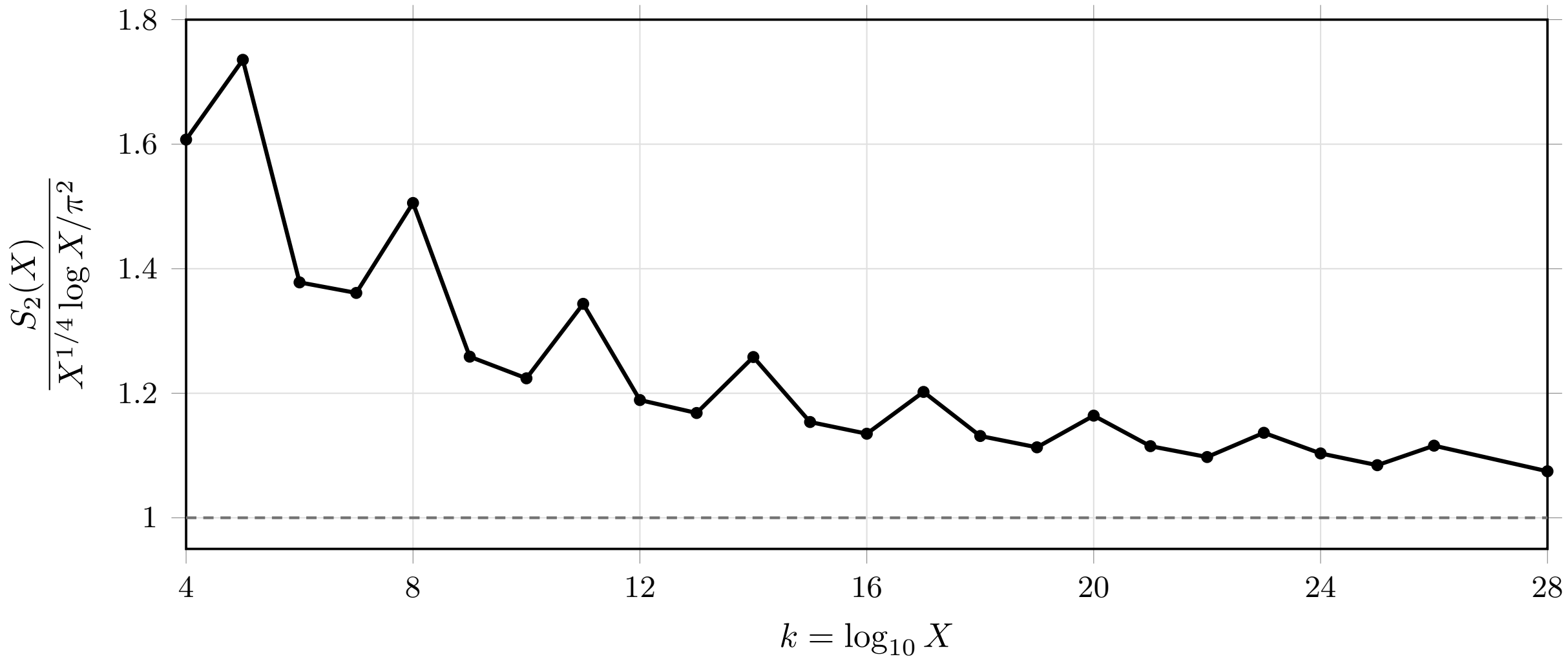
\begin{figure}[ht!]
\centering
\begin{tikzpicture}
  \begin{axis}[
    epsteinaxis,
    ymin=0.95, ymax=1.80,
    ytick={1.0,1.2,1.4,1.6,1.8},
    ylabel={$\displaystyle\frac{S_2(X)}{X^{1/4}\log X/\pi^2}$},
  ]
    \addplot[black!55,densely dashed,thick,domain=4:28]{1};
    \addplot[black,very thick,mark=*,mark size=1.15pt]
      table[x=k,y=ratio]{\steinData};
  \end{axis}
\end{tikzpicture}
\caption{Ratio of exact counts to the leading asymptotic. The dashed
horizontal line marks the limiting value $1$.}
\label{fig:epstein-ratio}
\credit{I input student-produced verified data to Astra, and Astra generated TikZ and PGFPlots code using said data. I manually adjusted the axis limits, tick marks, and plot styling for clarity and readability. Produced in 2026.}
\end{figure}

% \input{sections/discussion}
% \clearpage
\section{Remarks}

In Richard Guy's collection \cite{unsolved_problems} of unsolved problems in number theory, Problem E26 asks whether the losing positions or the winning positions of Epstein's game make up a positive-density subset of the positive integers. Our asymptotic settles this for the remoteness-2 positions: the E. P. Stein numbers form an infinite but density-zero subset of the positive integers, since

\[\lim_{X \to \infty} \frac{S_2(X)}{X}
  = \lim_{X \to \infty} \left(\frac{\log X}{\pi^2 X^{3/4}} + O(X^{-3/4})\right) = 0.\]

A related problem, A7(18), in Richard Nowakowski's \textit{Unsolved Problems in Combinatorial Games} \cite{polish} explicitly asks for asymptotic proportions of winning, losing, and draw positions. Although we have determined an asymptotic for remoteness-2 numbers, the general problem of asymptotic proportions for all positions remains open.

Further questions to investigate include whether $O(X^{1/4})$ is the best error obtainable and whether useful arithmetic characterizations exist for positions of higher remoteness in Epstein's game.

\clearpage
\bibliographystyle{ieeetr}
\bibliography{sections/references}

\end{document}